# Non-expansive maps between unit balls that send the sphere onto the sphere are linear isometries

Kaarel August Kurik
University of Maryland, College Park, MD 20742, USA
kaarelaugustkurik@gmail.com

**Abstract.** Let $X, Y$ be real Banach spaces with closed unit balls $B_X, B_Y$ and unit spheres $S_X, S_Y$. We prove that every non-expansive map $F : B_X \to B_Y$ with $F(S_X) \subseteq S_Y \subseteq F(B_X)$ is the restriction of a surjective linear isometry $X \to Y$. Consequently, a non-expansive bijection of $B_X$ onto itself that fails to be an isometry must map some point of $S_X$ into the open unit ball, so the plasticity of $B_X$ reduces to showing that non-expansive bijections preserve norm 1. The proof shows that $F$ is Fréchet differentiable at 0 with a surjective linear isometry $T$ as its derivative, and then proves $T^{-1}F = \mathrm{id}$ using a quantitative rigidity lemma for non-expansive self-maps $H$ of $B_X$ satisfying $\|a + H(a)\| = 2$ for $a \in S_X$; the latter relies on the Baillon–Bruck rate of asymptotic regularity for Krasnosel'skiĭ–Mann iterations. As a corollary, the identity is extremal in the set of non-expansive self-maps of $B_X$ for every Banach space $X$, which together with a result of Bargetz, Dymond and Pirk shows that all surjective isometries of $B_X$ are extremal, answering a question of theirs, including the case of $L^1[0, 1]$.

## 1 Introduction

The notion of a *plastic metric space* was introduced in [1] to distinguish those metric spaces $M$ for which every non-expansive (i.e. 1-Lipschitz) bijection of signature $M \to M$ is an isometry; every totally bounded metric space is plastic [1], [2]. A line of study into the plasticity of the closed unit balls of Banach spaces was initiated in [3], where it was shown that the closed unit ball of every strictly convex Banach space is plastic. A number of additional classes of Banach space have since been shown to have plastic unit balls, among them $\ell_1$ [4], $\ell_1$-sums of strictly convex spaces [5], spaces whose unit sphere is the union of its finite-dimensional polyhedral extreme subsets [6], $\ell_\infty$-sums of two strictly convex spaces and $\ell_1 \oplus_2 \mathbb{R}$ [7], $c$ [8], $C(K)$ for compact metrizable $K$ with finitely many accumulation points [9], and $\ell_\infty$ [10]; see also [11], [12] for non-expansive bijections between the unit balls of different spaces. It is currently unknown whether the closed unit ball of every Banach space is plastic.

The main result of this article is that when $X, Y$ are Banach spaces, $F : B_X \to B_Y$ is a non-expansive map and $F(S_X) \subseteq S_Y \subseteq F(B_X)$, then $F$ is the restriction of a surjective linear isometry of signature $X \to Y$. In the case $Y = X$, this is a substantial strengthening of Lemma 2.5 from [3], dropping the requirements of bijectivity, homogeneity and oddness for $F$. The new result tells us that any counterexample to the plasticity conjecture for closed Banach balls must strictly decrease norm at some input.

## 2 Statements of results

**Theorem 1.** *Let $X, Y$ be real Banach spaces with closed unit balls $B_X, B_Y$ and unit spheres $S_X, S_Y$. Let $F : B_X \to B_Y$ be a 1-Lipschitz map with $F(S_X) \subseteq S_Y$ and $S_Y \subseteq F(B_X)$. Then $F$ is the restriction of a surjective linear isometry $X \to Y$.*

**Lemma 2.** *Let $X$ be a real Banach space and let $H : B_X \to B_X$ satisfy*

$$\|a + H(a)\| = 2 \qquad (a \in S_X).$$

*If $H$ is $(1+\varepsilon)$-Lipschitz with $\varepsilon \geq 0$, then*

$$\sup_{x \in B_X} \|H(x) - x\| \leq \min\left(2, 4\varepsilon^{\frac{1}{3}}\right).$$

*In particular, if $H$ is non-expansive, then $H$ is the identity.*

## 3 Proof of Lemma 2

For $\varepsilon = 0$, the conclusion can also be derived from Theorem 3.3 in [13] (see also Theorem 9.4 in [14]); the argument below gives a quantitative version.

We may assume $X \neq \{0\}$, the statement being trivial otherwise.

### 3.1 Stability estimate

We apply the following theorem regarding Krasnosel'skiĭ–Mann iterates.

**Theorem** (Cominetti–Soto–Vaisman [15, Theorem 1])**.** *Let $C$ be a bounded convex subset of a normed space, let $K : C \to C$ be non-expansive and let $\alpha \in (0, 1)$. Then the iterates $x_n = (1-\alpha)x_{n-1} + \alpha K(x_{n-1})$ starting at any $x_0 \in C$ satisfy*

$$\|x_n - K(x_n)\| \leq \frac{\operatorname{diam} C}{\sqrt{\pi n \alpha(1-\alpha)}}.$$

Fix $\alpha = \frac{1}{8}$. Let $K : B_X \to B_X$ be nonexpansive and put

$$P = (1-\alpha)I + \alpha K, \qquad d_0 = \sup_{a \in S_X} (1 - \|P(a)\|).$$

Define the radial extension

$$\hat{K}(x) = r_x K(x/r_x), \qquad r_x = \max(1, \|x\|).$$

Note that $\hat{K}$ is 3-Lipschitz. With $x, y \in X$ and $r_x \geq r_y$:

$$\begin{aligned}
\left\|\hat{K}(x) - \hat{K}(y)\right\| &= \left\|r_x K(x/r_x) - r_y K\big(y/r_y\big)\right\| \\
&\leq \left(\left\|r_x K(x/r_x) - r_x K\big(y/r_y\big)\right\| + \left\|r_x K\big(y/r_y\big) - r_y K\big(y/r_y\big)\right\|\right) \\
&\leq r_x \left\|x/r_x - y/r_y\right\| + \left|r_x - r_y\right| \\
&= \left\|x - \frac{r_x}{r_y} y\right\| + \left|r_x - r_y\right| \\
&\leq \|x - y\| + \|y\| \left|1 - \frac{r_x}{r_y}\right| + \left|r_x - r_y\right| \\
&\leq \|x - y\| + 2\left|r_x - r_y\right|
\end{aligned}$$

$$\leq 3\|x-y\|.$$

Hence $\alpha\left(\hat{K}-I\right)$ is $\frac{1}{2}$-Lipschitz, from which $\hat{P} := I+\alpha\left(\hat{K}-I\right)$ is a bijection on $X$ by the Banach fixed point theorem. Note that $\hat{P}$ is the radial extension of $P$.

For $\|x\|>1$, we have

$$\left\|\hat{P}(x)\right\| = \|x\|\|P(x/\|x\|)\| \geq \|x\|(1-d_0) > 1-d_0.$$

That is, $\hat{P}(X\setminus B_X) \subseteq X\setminus(1-d_0)B_X$, hence $(1-d_0)B_X \subseteq \hat{P}(B_X) = P(B_X)$.

Every point in $B_X$ is thus within $d_0$ of $P(B_X)$. Since $P$ is nonexpansive, induction gives

$$\sup_{z\in B_X} \operatorname{dist}(z, P^n(B_X)) \leq nd_0.$$

Now, the constant-step asymptotic regularity bound (due to Baillon and Bruck [16]; see the cited theorem 1 above) yields

$$\sup_{x\in B_X} \|P^n(x)-K(P^n(x))\| \leq \frac{2}{\sqrt{\pi n\alpha(1-\alpha)}}.$$

Finally, with $\delta>0$ arbitrary, approximate an arbitrary $z\in B_X$ by a point $x\in P^n(B_X)$ such that $\|x-z\|\leq nd_0+\frac{\delta}{2}$ and obtain

$$\begin{aligned}\|K(z)-z\| &\leq \|K(z)-K(x)\|+\|K(x)-x\|+\|x-z\| \\ &\leq 2\|x-z\|+\frac{2}{\sqrt{\pi n\alpha(1-\alpha)}} \\ &\leq (2nd_0+\delta)+\frac{2}{\sqrt{\pi n\alpha(1-\alpha)}},\end{aligned}$$

from which

$$\sup_{z\in B_X}\|K(z)-z\| \leq 2nd_0+\frac{2}{\sqrt{\pi n\alpha(1-\alpha)}}.$$

### 3.2 Using the stability estimate

Set

$$L=1+\varepsilon, \qquad K=H/L, \qquad d=1-L^{-1}.$$

The map $K$ is nonexpansive, so the analysis from the previous section applies. With $\alpha, P$ as above and $a\in S_X$, the condition $\|a+H(a)\|=2$ implies that $a, H(a)$ are simultaneously normed by some norm-1 functional, which can be used to obtain

$$\|P(a)\| = 1-\alpha d \qquad (a\in S_X).$$

Thus $d_0=\alpha d$. Applying the stability estimate along with $\|H-K\|_\infty \leq d$, we obtain

$$\|H-I\|_\infty \leq d(1+2n\alpha)+\frac{2}{\sqrt{\pi n\alpha(1-\alpha)}} = \left(\frac{n}{4}+1\right)\frac{\varepsilon}{1+\varepsilon}+\frac{16}{\sqrt{7\pi n}}.$$

If $\varepsilon=0$, we let $n\to\infty$ to obtain $H=I$. If $0<\varepsilon\leq 1$, we choose $n=\left\lceil 4\varepsilon^{-\frac{2}{3}}\right\rceil$. We bound

$$\left(\frac{n}{4}+1\right)\frac{\varepsilon}{1+\varepsilon} \leq \left(\frac{n}{4}+1\right)\varepsilon \leq \frac{4\varepsilon^{-\frac{2}{3}}+5}{4}\varepsilon \leq \frac{5}{4}\varepsilon + \varepsilon^{\frac{1}{3}}$$

and

$$\frac{16}{\sqrt{7\pi n}} \leq \frac{16}{\sqrt{7\pi \cdot 4\varepsilon^{-\frac{2}{3}}}} = \frac{8}{\sqrt{7\pi}}\varepsilon^{\frac{1}{3}} < \frac{7}{4}\varepsilon^{\frac{1}{3}},$$

hence

$$\|H-I\|_\infty \leq \frac{11}{4}\varepsilon^{\frac{1}{3}} + \frac{5}{4}\varepsilon \leq 4\varepsilon^{\frac{1}{3}}.$$

For $\varepsilon > 1$, just use the trivial bound $\|H-I\|_\infty \leq \|H\|_\infty + \|I\|_\infty \leq 2$. □

## 4 Proof of Theorem 1

In the following, we assume $X, Y$ are nonzero spaces. This is no loss of generality: if $Y = \{0\}$, then $S_Y = \varnothing$, so $F(S_X) \subseteq S_Y$ forces $S_X = \varnothing$, i.e. $X = \{0\}$. If $X = \{0\}$, then $S_Y \subseteq F(B_X) = \{F(0)\}$, and since $S_Y = -S_Y$ cannot be a singleton, $Y = \{0\}$. In either case the statement is trivial.

### 4.1 Norm preservation, defect estimates

Let $F$ be as in the statement of Theorem 1.

We have

$$1 + \|F0\| = \sup_{y \in S_Y} \|F0 - y\| \leq \sup_{x \in B_X} \|F0 - Fx\| \leq 1,$$

hence $F0 = 0$ and $\|Fx\| \leq \|x\|$.

For $v \in S_X$ and $r \in [0,1]$, we have

$$\|F(rv)\| \geq \|Fv\| - \|Fv - F(rv)\| \geq 1 - \|v - rv\| = r \geq \|F(rv)\|,$$

so

$$\|Fx\| = \|x\| \quad (x \in B_X).$$

Next, we're going to be working with vectors formed as finite sums $D = \sum_i c_i F(x_i)$ with $c_i > 0$. Later on we'll choose $c_i, x_i$ such that if $F$ were linear, $D$ would equal zero exactly: these will be expressions like $uF(tv) + tF(-uv)$, or $F(tv) + F(tw) + F(-t(v+w))$. Consequently $D$ stands for "defect". For our first estimates, we don't need this detail.

First, choose $a \in S_X$ such that $\|D\|Fa = -D$. This is possible: if $D = 0$, any $a \in S_X$ will do, and otherwise $-D/\|D\| \in S_Y \subseteq F(B_X)$, and any preimage of it lies in $S_X$ by norm preservation. If $b \in B_X$ and $\|a - b\| \leq \frac{1}{4}$, choose any norm-one functional $\psi$ norming $Fb$, which is possible since $\|Fb\| = \|b\| \geq \frac{3}{4}$. By norm preservation for $F$, we have

$$\psi(Fa) \geq \psi(Fb) - \|Fa - Fb\| = \|Fb\| - \|Fa - Fb\| \geq \|b\| - \|a - b\| \geq \frac{3}{4} - \frac{1}{4} = \frac{1}{2}.$$

For every $i$, comparison with the same $b$ gives

$$\psi(Fx_i) \geq \|Fb\| - \|Fb - Fx_i\| \geq \|b\| - \|b - x_i\|.$$

Summing these bounds weighted by $c_i$, we obtain

$$\psi(D) \geq \left(\sum_i c_i\right)\|b\| - \sum_i c_i\|b - x_i\|.$$

We have $\psi(D) = -\|D\|\psi(Fa) \leq -\frac{\|D\|}{2}$, hence

$$\frac{\|D\|}{2} \leq \sum_i c_i\|b - x_i\| - \left(\sum_i c_i\right)\|b\|.$$

This bound holds for all $b \in B_X$ with $\|b - a\| \leq \frac{1}{4}$.

**4.2 Obtaining a directional derivative at zero**

Fixing $v \in S_X$ and $0 < t, u \leq 1$, we first aim to prove

$$\left\|\frac{F(tv)}{t} + \frac{F(-uv)}{u}\right\| \leq 8(t + u). \tag{1}$$

Set the defect $D = uF(tv) + tF(-uv)$. If $D = 0$, we're already done. Otherwise, choose $a$ as in the previous section, set $h = \frac{1}{8}$, and consider the family

$$b_s = \frac{7}{8}a + sv \qquad (-h \leq s \leq h).$$

These belong to $B_X$ and satisfy $\|a - b_s\| \leq \frac{1}{4}$. Define

$$g(s) := \left\|\frac{7}{8}a + sv\right\|.$$

This function is convex and 1-Lipschitz.

We have $\|b_s - tv\| = g(s - t)$, $\|b_s - (-uv)\| = g(s + u)$ and $\|b_s\| = g(s)$. Thus the bound on the defect norm obtained in the previous section becomes

$$\frac{\|D\|}{2} \leq ug(s - t) + tg(s + u) - (t + u)g(s).$$

Define the interval average of $g$ as

$$G(r) := \frac{1}{2h}\int_{-h}^{h} g(s + r)\,\mathrm{d}s.$$

We have

$$G'(r) = \frac{g(r + h) - g(r - h)}{2h}, \qquad |G'(r) - G'(q)| \leq \frac{|r - q|}{h}.$$

Hence, for every real $r$, we have

$$G(r) \leq G(0) + rG'(0) + \frac{r^2}{2h}.$$

We can average the RHS of our bound on $\frac{\|D\|}{2}$ over $s \in [-h, h]$ and apply the above inequality at $-t, u$ to obtain

$$\frac{\|D\|}{2} \leq uG(-t) + tG(u) - (t + u)G(0) \leq \frac{ut^2 + tu^2}{2h}.$$

We may divide this by $tu$ and use $h = \frac{1}{8}$ to obtain (1).

Let's write $W_t(v) := \frac{F(tv)}{t}$. With $u = t$, (1) gives

$$\|W_t(v) + W_t(-v)\| \leq 16t.$$

Thus we obtain

$$\|W_t(v) - W_u(v)\| \leq \|W_t(v) + W_u(-v)\| + \|W_u(-v) + W_u(v)\| \leq 8t + 24u.$$

Since $Y$ is complete, this tells us $\lim_{t\downarrow 0} W_t(v)$ exists. We note that for $k > 0$, $W_{kt}(v) = \frac{1}{k} W_t(kv)$, so the limit exists for any $v \in X$. We define

$$Tv := \lim_{t\downarrow 0} W_t(v) \qquad (v \in X)$$

and note that $T$ is positively homogeneous: $T(kv) = \lim_{t\downarrow 0} W_t(kv) = \lim_{t\downarrow 0} kW_{kt}(v) = kT(v)$.

By the bound on $\|W_t(v) + W_t(-v)\|$, we have $Tv + T(-v) = 0$. By passing to the limit of $\|W_t(v)\|$ and $\|W_t(v) - W_t(w)\|$, we obtain

$$\|Tv\| = \|v\|, \qquad \|Tv - Tw\| \leq \|v - w\|.$$

Letting $u \downarrow 0$ in (1), we have $\left\| \frac{F(tv)}{t} - Tv \right\| \leq 8t = 8\|tv\|$, or equivalently,

$$\|Fx - Tx\| \leq 8\|x\|^2 \qquad (x \in B_X) \tag{2}$$

by choosing $v := \frac{x}{\|x\|}$.

### 4.3 Obtaining additivity by averaging over two directions

Fix $v, w \in B_X$ and $0 < t \leq \frac{1}{2}$. Set

$$D_t := F(tv) + F(tw) + F(-t(v + w)).$$

If $D_t \neq 0$, choose $a$ with $Fa = -D_t/\|D_t\|$. Set $h = \frac{1}{16}$ and

$$g(s, r) := \left\| \frac{7}{8}a + sv + rw \right\|.$$

For $|s|, |r| \leq h$, the associated points belong to $B_X$ and are within $\frac{1}{4}$ of $a$. The defect bound gives

$$\frac{\|D_t\|}{2} \leq g(s - t, r) + g(s, r - t) + g(s + t, r + t) - 3g(s, r). \tag{3}$$

The function $g$ is 1-Lipschitz for the $\ell_1$ norm on $\mathbb{R}^2$. Define the average

$$G(z) := \frac{1}{4h^2} \int_{-h}^{h} \int_{-h}^{h} g(z_1 + s, z_2 + r) \, \mathrm{d}s \, \mathrm{d}r.$$

The first partials of $G$ are $\frac{1}{h}$-Lipschitz for the $\ell_1$ norm. Integration along the segment from 0 to $z$ therefore gives

$$G(z) \leq G(0) + \nabla G(0) \cdot z + \frac{\|z\|_1^2}{2h}. \tag{4}$$

Averaging (3) and applying (4) at $(-t, 0), (0, -t), (t, t)$, we obtain

$$\frac{\|D_t\|}{2} \le \frac{t^2 + t^2 + (2t)^2}{2h} = 48t^2.$$

It suffices to divide by $t$ and use the norm limits from the previous section to conclude

$$Tv + Tw + T(-v - w) = 0.$$

Given that $T$ is odd and positively homogeneous, this suffices to show that $T$ is linear, and since it preserves norms, it's an isometric embedding. Inequality (2) tells us that $T$ is the Fréchet derivative of $F$ at 0.

### 4.4 Showing that $T$ is a surjection

Let's first show that the norming functionals of $Tx$ include all the norming functionals of $Fx$.

Let $a \in S_X$ and $\psi \in S_{Y^*}$ such that $\psi(Fa) = 1$. For $0 < t \le 1$, we have

$$t \ge \psi(F(ta)) \ge \psi(Fa) - \|Fa - F(ta)\| \ge 1 - (1 - t) = t.$$

Thus

$$\psi(Ta) = \lim_{t \downarrow 0} \psi\left(\frac{F(ta)}{t}\right) = 1. \tag{5}$$

Next, suppose $T(X) \ne Y$. Since $T(X)$ is closed, Hahn-Banach gives some $\varphi \in S_{Y^*}$ vanishing on $T(X)$. By Bishop-Phelps [17], norm-attaining functionals are dense in $Y^*$, so we may choose some $\psi$ within $\frac{1}{3}$ of $\varphi$ and consider its normalization $\hat{\psi}$. We have $\left\|\hat{\psi} - \varphi\right\| \le \frac{2}{3}$, and $\hat{\psi}$ attains its norm at some $y \in S_Y$. By sphere coverage, we have $y = Fa$ with $a \in S_X$, and (5) yields $\hat{\psi}(Ta) = 1$. However, $\varphi(Ta) = 0$, so $1 \le \left\|\hat{\psi} - \varphi\right\| \le \frac{2}{3}$, a contradiction. Hence $T(X) = Y$.

### 4.5 Common face rigidity

The map $H = T^{-1}F : B_X \to B_X$ is nonexpansive. For $a \in S_X$, apply (5) with $\psi$ norming $Fa$. Then $\|Ta + Fa\| = 2$, hence $\|a + Ha\| = 2$.

From Lemma 2, $H = I$, hence $F = T$ on $B_X$. □

## Use of AI

The human author of this article attributes the preceding results entirely to the work of OpenAI's and Anthropic's language models (primarily GPT-6 Astra and Claude Opus 5.5), and vouches for their correctness and significance. The abstract, bibliography and typesetting were prepared with the assistance of Claude Opus 5.5 and checked by the author.

## Formalization

Theorem 1 and Lemma 2 have both been proved formally in the Lean proof assistant, with their statements and proofs hosted on the Prove2Me platform, at

- https://prove2.me/theorems/7db30a91-7e0a-4b33-b52b-5e58b046b84f (Theorem 1, with the benign restriction that $X, Y$ are nontrivial),
- https://prove2.me/theorems/c93402b5-174f-4a0d-921d-6576eca38227 (Lemma 2)

## Appendix A. Extremality of the identity

Lemma 2 yields the following corollary:

**Corollary 3.** *Let $X$ be an arbitrary real Banach space. The identity is extremal in the space of non-expansive mappings $B_X \to B_X$.*

*Proof.* Suppose there exist non-expansive maps $f, g : B_X \to B_X$ such that

$$\mathrm{id} = \lambda f + (1-\lambda) g, \quad \lambda \in (0,1).$$

Fix $a \in S_X$ and let $\varphi$ be a norm-one functional such that $\varphi(a) = 1$. Then

$$1 = \varphi(a) = \lambda\varphi(f(a)) + (1-\lambda)\varphi(g(a)) \leq \lambda\varphi(f(a)) + (1-\lambda) \leq 1,$$

hence $\varphi(f(a)) = 1$, so $\|a + f(a)\| = 2$. Lemma 2 then yields $f = \mathrm{id}$, and consequently $g = \mathrm{id}$ as well. □

In conjunction with Lemma 3.1 in [18], Corollary 3 implies that every surjective isometry $B_X \to B_X$ is extremal in the space of non-expansive maps $B_X \to B_X$. This is a generalization of Theorem 1.1 from the same paper and answers the question raised there of whether surjective isometries are extremal, including the case of $L^1[0,1]$, which the authors left open.